\documentclass{article}
\usepackage[labelsep=none]{caption}
\usepackage{amssymb,color,graphicx,subfig}
\usepackage[utf8]{inputenc}
\newcommand\N{\mathbb{N}}
\newcommand\Z{\mathbb{Z}}
\newcommand\R{\mathbb{R}}
\newcommand\T{\mathbb{T}}
\newcommand\La{\Lambda}
\newcommand\M{\mathcal{M}}
\newcommand\vol{\mathrm{vol}}
\newcommand\vect[1]{{\mathbf #1}}
\newcommand\marca[1]{}
\newtheorem{prop}{Proposition}
\newtheorem{lemma}[prop]{Lemma}
\newtheorem{theo}[prop]{Theorem}
\newtheorem{cor}[prop]{Corollary}
\newenvironment{demo}{\noindent{\bf Proof.}
}{\hfill\(\blacksquare\)\\ }
\usepackage[a4paper,includehead,footskip=1cm,footnotesep=1cm,marginparwidth=1cm,top=15mm,bottom=2cm,left=20mm,right=20mm]{geometry}
\usepackage[colorlinks=true,filecolor=blue,citecolor=blue,linkcolor=blue,urlcolor=blue]{hyperref}
\usepackage[sort&compress,square]{natbib}
\begin{document}
\begin{center}
\parbox{11cm}{\begin{center}\Large\bf Full support of the Kasteleyn operator associated with a bipartite toroidal graph\end{center}}\vspace{3mm}

\renewcommand{\thefootnote}{\fnsymbol{footnote}}

{\large Álvar Ibeas Martín\footnote{alvar.ibeas@unican.es}}\\[3mm]
{\large Universität Wien}
\end{center}

\noindent{\bf Abstract.} A perfect matching in a bipartite graph
embedded on a torus defines a height function on the graph's faces and
an associated height change vector in \(\Z^2\). These matchings are
enumerated by a combination of four evaluations of a bivariate
Laurent polynomial, called Kasteleyn operator, whose coefficient of
bidegree \((i,j)\) is, up to the sign, the number of perfect matchings
with height change \((i,j)\). Therefore the Newton polygon of the
Kasteleyn operator is the convex hull of the height change vectors. In
this article, we prove that any point with integer coordinates in that
polygon is realized by a perfect matching.

\section{Introduction}

\begin{figure}[h!]
\begin{center}
\subfloat[]{\label{fig_uv}\includegraphics[height=5cm]{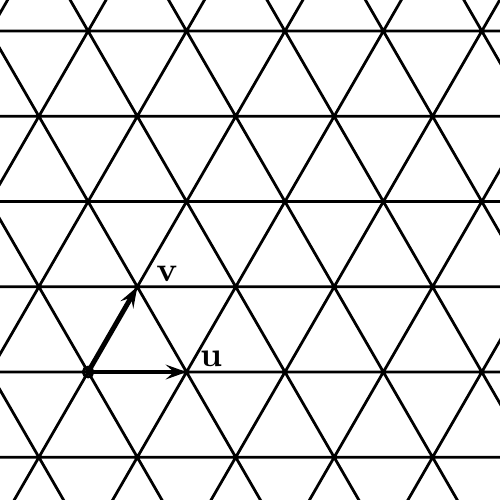}}\hspace{5mm}
\subfloat[]{\label{fig_plt}\includegraphics[height=5cm]{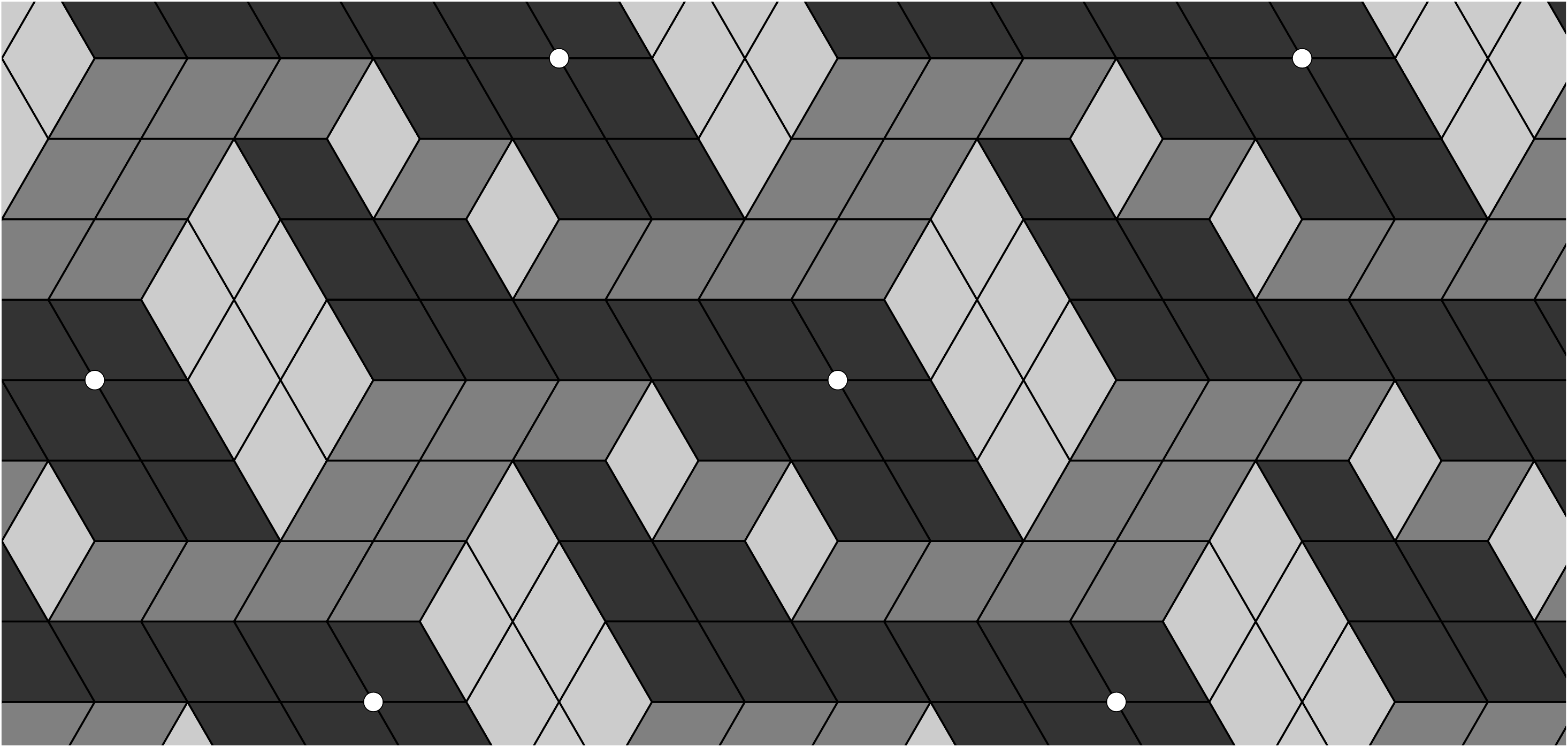}}
\caption{}
\end{center}
\end{figure}

Consider the set of doubly periodic lozenge tilings of the plane, for
a given period \(\mathcal{L}=(\vect u|\vect v)B\Z^2\), where \(B\) is
an integer \(2\times 2\) matrix with nonzero determinant and \(\vect
u,\vect v\) are the (column) vectors depicted in
Figure~\ref{fig_uv}. For instance, the lozenge tiling shown in
Figure~\ref{fig_plt} has the period \(\Z\langle 3\vect u+4\vect
v,-5\vect u+4\vect v\rangle\).

The natural identification of lozenge tilings with stacks of cubes
induces a height function that labels the points of
\(\mathcal{L}_0:=\Z\langle\vect u,\vect v\rangle\). Being the
considered tilings invariant through translations by vectors of the
sublattice \(\mathcal{L}\), such a translation produces a constant
increment in the height function. For example, in
Figure~\ref{fig_plt}, the vector \(3\vect u+4\vect v\) reduces the
height by 2 and the vector \(-5\vect u+4\vect v\) keeps it unchanged,
so we can associate the pair \((-2,0)\) to this tiling. Not all the
tilings with this period produce the same height change. For
instance, the three constant tilings consisting of lozenges of the
same type (\rotatebox[origin=c]{-60}{\(\lozenge\)},
\rotatebox[origin=c]{58}{\(\lozenge\)}, and \(\lozenge\),
respectively) give the following height changes: \((1,9)\),
\((-11,-3)\), and \((10,-6)\).
\begin{figure}[h!]
\begin{center}
\includegraphics[height=4cm]{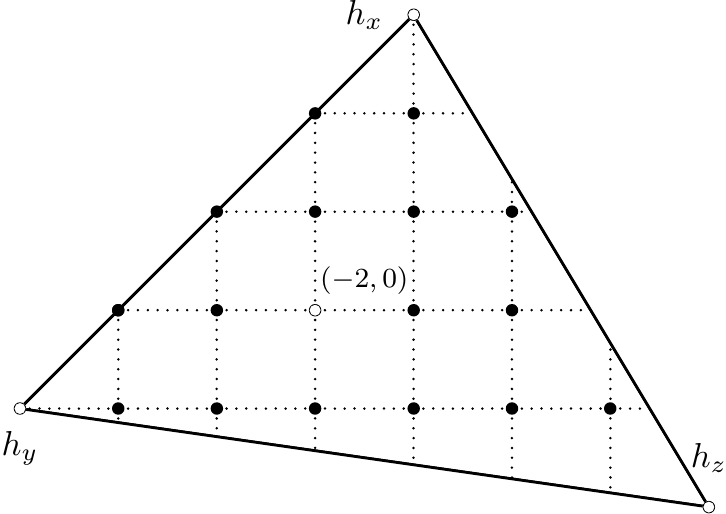}
\end{center}
\end{figure}

For a general period \(\mathcal{L}\subseteq\mathcal{L}_0\) (and a
fixed basis), the possible height change vectors can be
determined. Let \(h_x,h_y,h_z\) be the vectors associated to the
constant tilings. If, for a \(\mathcal{L}\)-periodic tiling, the
amount of lozenges of each type in a fundamental cell of the period is
\((x,y,z)\), its height change is the following convex combination:
\(\frac{\vol(\mathcal{L}_0)}{\vol(\mathcal{L})}(xh_x+yh_y+zh_z)=\frac{1}{|\det(B)|}(xh_x+yh_y+zh_z)\),
which is a point in the triangle determined by \(h_x\), \(h_y\), and
\(h_z\). Moreover, its difference to any of these vertices lies in
\(3\Z^2\), for the height difference between two points of
\(\mathcal{L}_0\) is determined modulo 3, independently of the tiling
or the period.

Conversely, for any point \(h\) in the triangle with coordinates in
\(h_x+3\Z^2\), there exists an \(\mathcal{L}\)-periodic tiling whose
height change vector is \(h\). The goal of this article is to give
a proof of this fact in a more general setting, which considers
perfect matchings in a general bipartite graph embedded on a torus
instead of periodic lozenge tilings.

A slight modification of the construction explained above is a
particular case of the height function defined by
\citet{Ken_Oko_She_06} in the general case (see also
\citep{Thu_90}). Summing up, a pair of perfect matchings
\((\omega_0,\omega)\) defines a set of circuits (called
\emph{transition graph}) on the torus whose homology type
\(h_{\omega_0}(\omega)\in\Z^2\) determines the height change
vector. Details follow in Subsection~\ref{sec_hei}.

For a given bipartite toroidal graph and a fixed base perfect matching
\(\omega_0\), the associated \emph{Newton polygon} is the convex hull
of the set of homology types \(h_{\omega_0}(\omega)\) of the
transition graphs formed by \(\omega_0\). It receives this name
because it coincides with (a translate of) the Newton polygon of a
bivariate Laurent polynomial \(P(w,z)\), namely, the determinant of a
weighted Kasteleyn-Percus matrix (see \citep{Ken_Oko_She_06}). The
absolute value of a coefficient of \(P\) is the number of perfect
matchings with a corresponding height change. In
Section~\ref{sec_main} we prove that any point with integer
coordinates in the Newton polygon is realized by some perfect matching,
so that the support (i.e. the set of occurring monomials) of the
Kasteleyn operator is maximal. \marca{\citep{Gon_Ken_11}}

\section{Preparations}

In this section we collect some basic facts that will be needed in the
proof of our theorem. Firstly, we explain the formalization of the
height functions described in the introduction, following
\citep{Ken_L}. Then, we discuss the homology of a set of knots in the
torus. Finally, we relate the existence of circuits in circulant
digraphs to the visibility of lattice points.

\subsection{Height functions}\label{sec_hei}

In order to fix the class of graphs whose perfect matchings are
assigned height functions, let us recall some facts about topological
graphs. We refer to \citep{Gro_Tuc_87,Moh_88} for a detailed
treatment. A graph \(G=(V,E)\) (possibly with loops or multiple edges)
is endowed with a natural topology. An \emph{embedding} of a graph
\(G\) on a surface \(S\) is a mapping
\(i:G\lhook\joinrel\relbar\joinrel\rightarrow S\) such that the
restriction \(i:G\lhook\joinrel\relbar\joinrel\rightarrow i(G)\) is a
homeomorphism (with \(i(G)\) endowed with the subspace topology). We
usually identify a graph and its image through an embedding, so that
\(G\subset S\). Note that, as we deal with infinite graphs,
the usual definition of embedding as an injective and continuous
mapping is not equivalent. An embedding is \emph{cellular} if the
complement of the graph in the surface is homeomorphic to a disjoint
union of open discs.

We use the term \emph{periodic graph} for a cellular embedding of a
graph on the plane that is periodic through integer-valued vectors,
i.e. the translation by any vector in \(\Z^2\) is a graph
automorphism. The existence of a cellular embedding on the plane
implies that the graph is locally finite. We will consider bipartite
periodic graphs, but will not require them to be connected.

\begin{figure}[h!]
\begin{center}
\includegraphics[width=3cm]{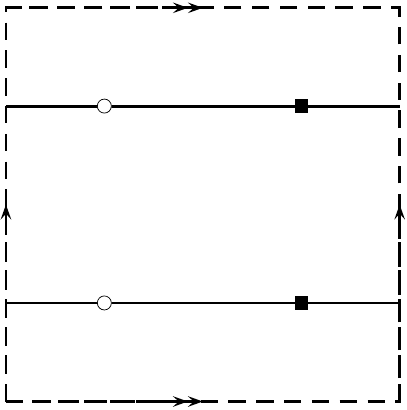}
\caption{}\label{ej_1}
\end{center}
\end{figure}

The projection of a periodic graph \(G\) on the torus \(\T^2=\R^2/\Z^2\) is
the embedding of a graph \(\hat G\). This graph is finite, for it is
locally finite in a compact surface. Note that the embedding need not
be cellular, as is shown in Figure~\ref{ej_1}.

We use the notation \(\M(G)\) for the set of \emph{perfect matchings}
or \emph{1-factors} of a graph \(G\), i.e.
\[\M(G):=\{\omega\subseteq E\ |\ \forall v\in V\ \exists! e\in\omega\
:\ v\in e\}.\]
In the following, we use the terms \emph{matching} and \emph{perfect
  matching} indistinctly, for we do not consider incomplete
matchings. We are interested in those matchings on \(G\) which are
compatible with the projection on the torus; or equivalently,
matchings on \(\hat G\).

For a bipartite periodic graph, we set a fixed orientation (from the
white to the black vertex) at every edge. Given two matchings
\(\omega,\omega'\in\M(\hat G)\), the \emph{transition graph}
\(\omega-\omega'\) is composed by the oriented edges of \(\omega\) and
the reversed edges of \(\omega'\). The connected components of a
transition graph are directed circuits (\emph{transition cycles},
according to \cite{Kas_63}) and pairs of vertices bidirectionally
linked, which we discard.  We assume that the matching set is nonempty
and fix a base matching \(\omega_0\in\M(\hat G)\). Let
\(\omega\in\M(\hat G)\) and let us define, by means of the transition
graph \(\omega-\omega_0\), a \emph{height function} on the faces of
\(G\), i.e. the vertices of the geometric dual \(G^*\). Note that
\(G^*\) is connected, although it need not be locally finite (and
therefore, it is not necessarily a cellular embedding on the
plane). Firstly, we choose a base face and assign the height 0 to
it. If there is an edge between two faces which does not occur in the
transition graph, both faces get the same height. If the dual edge
that goes from a face \(F_1\) to a face \(F_2\) is crossed from left
to right by the transition graph, we set \(h(F_2)=h(F_1)+1\).
\begin{figure}[h!]
\begin{center}
\includegraphics{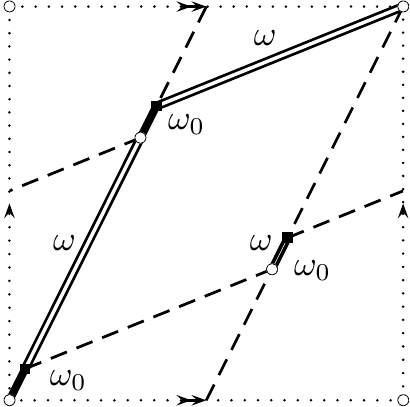}\hfill
\includegraphics{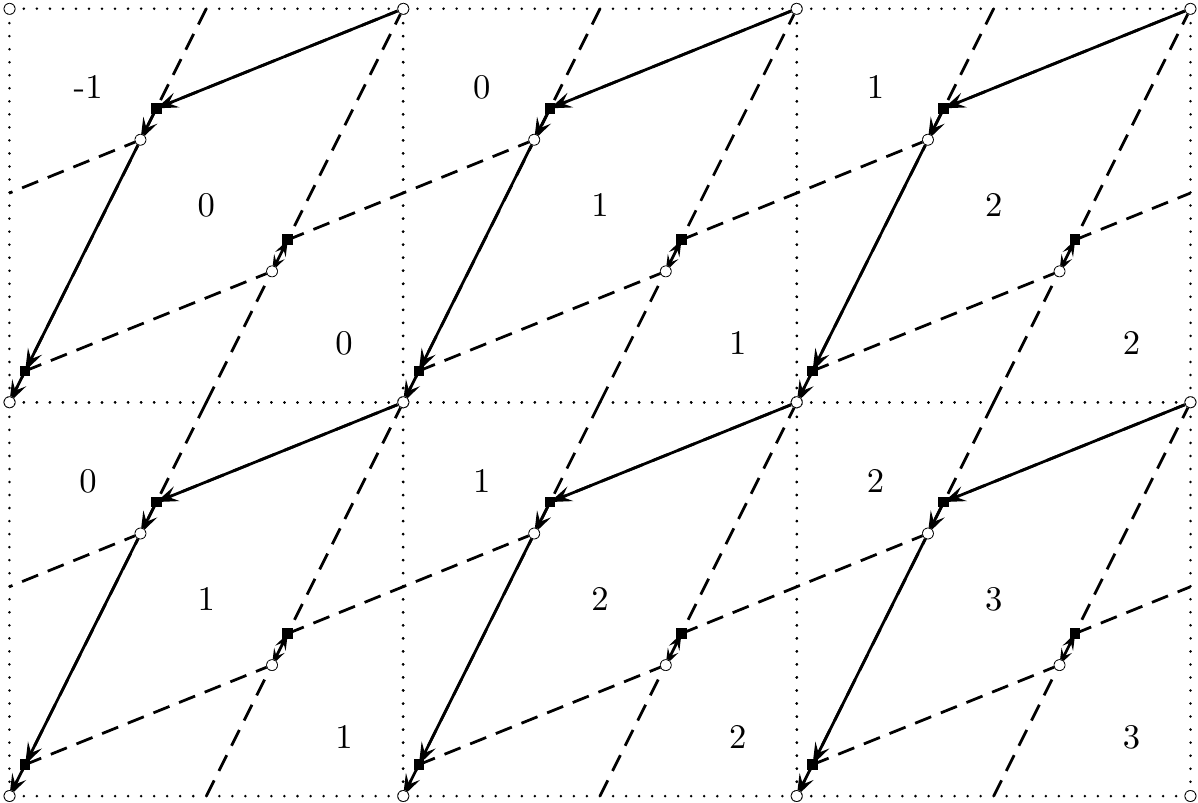}
\caption{}\label{ej_2}
\end{center}
\end{figure}

This process consistently defines an integer-valued function on the
faces of \(G\) (see an example in Figure~\ref{ej_2}). To see this,
consider the three free Abelian groups generated by the vertex set of
\(G\), the set of edges endowed with the orientation defined above,
and the set of (oriented according to the embedding) faces,
respectively. The sketch below represents the border homomorphisms
\(\partial_i\) and their transposed: the coborder homomorphisms
\(\delta_i\).
\begin{figure}[h!]
\begin{center}
\includegraphics{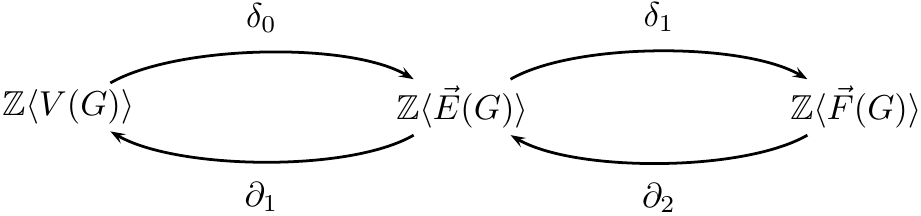}
\end{center}
\end{figure}

Writing \(\partial^*\) and \(\delta^*\) for the border and coborder
homomorphism (respectively) associated with the dual graph, we have
\(\partial_1^*=\delta_1\), \(\partial_2^*=\delta_0\). An element
\(\tau^*\in\Z\langle\vec E(G^*)\rangle\) defines the homomorphism
\[\begin{array}{cccc}\tau^*=\langle\tau,\cdot\rangle:&\Z\langle\vec
  E(G)\rangle&\longrightarrow&\Z\\
&\omega&\mapsto&\displaystyle\sum_{e\in\vec E}\tau(e)\omega(e),
\end{array}\]
and an analogue construction can be done for elements of
\(\Z\langle\vec F(G^*)\rangle\). With this notation, if \(F_1,F_2\in
F(G)\) and \(\tau^*\) is a walk from \(F_1\) to
\(F_2\) in the dual graph, we have
\(h(F_2)-h(F_1)=\langle\tau,\omega-\omega_0\rangle\). Now, if
\(\tau_1^*\) and \(\tau_2^*\) are walks in \(G^*\) with the same
endpoints, we have to see that
\(\langle\tau_1,\omega-\omega_0\rangle=\langle\tau_2,\omega-\omega_0\rangle\). Note
that \(\tau_1^*-\tau_2^*\) is a border (there is \(\rho^*\in\Z\langle\vec
F(G^*)\rangle\) such that \(\partial^*(\rho^*)=\tau_1^*-\tau_2^*\)) and
\(\omega-\omega_0\) is a cycle
\((\partial(\omega-\omega_0)=0)\). Therefore \(\langle\tau_1-\tau_2,\omega-\omega_0\rangle=\langle\delta(\rho),\omega-\omega_0\rangle=\langle\partial(\omega-\omega_0),\rho\rangle=0\).

On the other hand, if \(F\) and \(F'\) are faces of \(G\) and \(\vect
u\in\Z^2\), we have \(h(F+\vect u)-h(F)=h(F'+\vect u)-h(F')\). The
height function is determined, therefore, by its values on a system of
face representatives modulo \(\Z^2\) and the pair of height increments
\(\tilde h_{\omega_0}(\omega):=(h(F+(1,0))-h(F),h(F+(0,1))-h(F))\). As
before, it can be shown that the height change vector is determined
by the homology type of \(\omega-\omega_0\) in the torus. It is easy
to check that, if the homology type of \(\omega-\omega_0\) is
\(h_{\omega_0}(\omega)=(a,b)\in\Z^2\), we have \(\tilde
h_{\omega_0}(\omega)=(-b,a)\). For instance, in the example of
Figure~\ref{ej_2}, we have \(\tilde h_{\omega_0}(\omega)=(1,-1)\) and
\(h_{\omega_0}(\omega)=(-1,-1)\).

\subsection{Torus knots}\label{sec_top}

As we have seen, the height change of a matching (with respect to a
base matching) can be identified with the homology type of a set of
disjoint oriented copies of \(S^1\) on the torus. Let us collect some
results on torus knots (whose proofs can be found in~\citep{Rol_76})
for later use. Note that, unlike usually, we consider directed knots.
\begin{lemma}Let \(\vect u=(u_1,u_2)\in\Z^2\). There exists a torus
  knot \(c\) with homology type \(\vect u\) if and only if \(\vect u\)
  is a visible lattice point (i.e. \(\vect u=\vect 0\) or
  \(\gcd(u_1,u_2)=1\)).
\end{lemma}
Let \(f\in\mathrm{Aut}(\T^2)\) be a self-homeomorphism on the
torus. The functor \(\pi_1\) associates an automorphism of the
fundamental group \(\Z^2\) with it:
\[\mathrm{Aut}(\T^2){\stackrel{\pi_1}{\longrightarrow}}\mathrm{GL}(2,\Z).\]
This mapping is a group epimorphism, in particular, for every group
automorphism \(P\) on \(\Z^2\), there is a self-homeomorphism on the
torus whose effect on the homology types of the torus cycles is
determined by \(P\). Moreover:
\begin{lemma}\label{lem_b} Let \(c_1,c_2\) be two
  torus knots with nontrivial homology. Then there is a
  self-homeomorphism on the torus which maps \(c_1\) into \(c_2\).
\end{lemma}
\begin{lemma}Let \(c_1,c_2\) be two disjoint torus knots with
  nontrivial homology types. Then both homology types coincide or are
  opposite.
\end{lemma}
\begin{figure}[h!]
\begin{center}
\includegraphics[scale=.7]{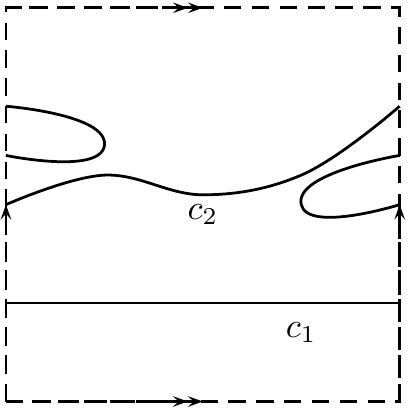}
\end{center}
\end{figure}
 \begin{demo}By Lemma~\ref{lem_b}, we can assume that the homology type
   of \(c_1\) is \((1,0)\). Therefore \(c_2\) is a knot on the
   cylinder \(\T^2\backslash c_1\), so its homology type must be
   \(\pm(1,0)\).
 \end{demo}
The following result is a direct consequence.

\begin{lemma}\label{lem_divide}Let \(G\) be a bipartite periodic graph and
  \(\omega_0,\omega\in\M(\hat G)\), such that \(\vect
  u=(u_1,u_2):=h_{\omega_0}(\omega)\in\Z^2\backslash\{\vect 0\}\). Set
  \(d:=\gcd(u_1,u_2)\). Then the transition graph consists of
  circuits with zero homology, \(P\) circuits with homology type
  \(\frac 1 d\vect u\), and \(N\) circuits with homology type
  \(\frac{-1}d\vect u\), with \(P-N=d\). In particular, if the
  transition graph consists of a single circuit, \(\vect u\) is a
  visible lattice point.
\end{lemma}

Note that if a transition graph \(\omega-\omega_0\) consists of
several circuits, the removal of some of them leads to another
transition graph \(\omega'-\omega_0\). Using this remark, we get:

\begin{lemma}\label{lem_divide_bis}
  Let \(G\) be a bipartite periodic graph and
  \(\omega_0,\omega\in\M(\hat G)\), such that \(\vect
  u=(u_1,u_2):=h_{\omega_0}(\omega)\in\Z^2\backslash\{\vect 0\}\). Set
  \(d:=\gcd(u_1,u_2)\). Then there is a matching \(\omega'\in\M(\hat
  G)\) such that the transition graph \(\omega'-\omega_0\) consists of
  one single circuit with homology type \(\frac 1 d\vect u\).
\end{lemma}

\subsection{Circulant digraphs}\label{sec_circ}

In the proof of the main result of this article we use a certain class
of toroidal graphs whose matchings can be identified with sets of
disjoint circuits in directed circulant graphs. Let us recall the
definition and some properties of these objects (see the survey by
\citet{Ber_Com_Hsu_95}).

For integers \(n,j_1,\ldots,j_r\) such that \(n\geq 1\), the
\emph{circulant graph} \(C(n;j_1,\ldots,j_t)\) is the Cayley graph
with vertex group \(\Z/n\Z\) and edges \(\{\{x,x+j_i\}\ |\ 1\leq i\leq
t\}\), where the \emph{jumps} \(j_i\) are considered as residue
classes modulo \(n\).
\captionsetup{labelsep=colon}
\begin{figure}[h!]
\begin{center}
\includegraphics[scale=.5]{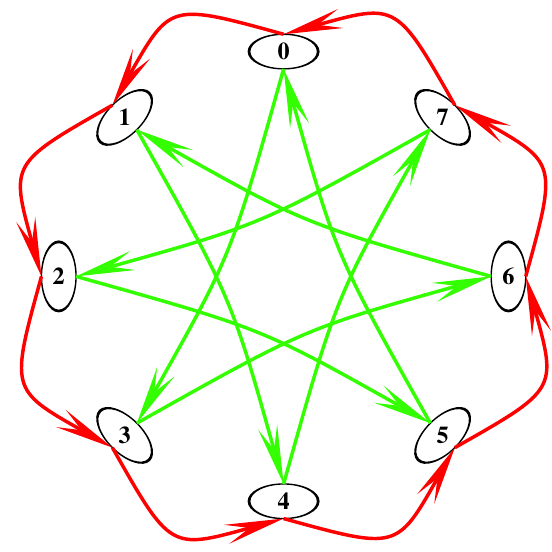}
\caption{\(\vec C(8;1,3)\)}
\end{center}
\end{figure}
\captionsetup{labelsep=none}

Note that we are including in the definition circulants with repeated
jumps (leading to multiple edges) and zero jumps (leading to
loops). We denote by \(\vec C(n;j_1,\ldots,j_t)\) the corresponding
directed graph. For instance, the digraph \(\vec C(n;1,1)\) consists
of a directed cycle of length \(n\) with duplicated arcs.

Let us consider a circulant digraph with two jumps: \(\vec
C(n;a,b)\). If a walk consists of \(u\) arcs of the type \((x,x+a)\)
and \(v\) arcs of the other type, we say that its \emph{abelianized}
is \((u,v)\). If \(x\in\Z/n\Z\) is the initial vertex of such a walk,
its terminal vertex is \(x+au+bv\). We can label, as is done in
\citep{Won_Cop_74}, a point \((u,v)\in\N^2\) by \(au+bv\): the
terminal vertex of a walk with abelianized \((u,v)\) and initial
vertex 0.

Note that the set of abelianized closed walks in a connected graph
\(\vec C(n;a,b)\) is the intersection of a 2-rank integer lattice with
volume \(n\) (which we call \emph{circuit lattice}) with
\(\N^2\). For instance, in the case \(\vec C(n;1,b)\), the circuit
lattice is \(\Z\langle(n,0),(-b,1)\rangle\).

We use the term \emph{lattice path} for a finite list of points in
\(\N^2\) such that the difference between two consecutive points is
\((1,0)\) or \((0,1)\). The labelling considered above allows the
identification of walks in the digraph with lattice paths. When there
are no two distinct points in a lattice path which are congruent
modulo the circuit lattice, the associated walk in the digraph is a
path (i.e. no vertex is visited twice). When the only pair of
congruent points modulo the circuit lattice consists of the endpoints
of the lattice path, the associated walk is a circuit.

\begin{lemma}\label{lema_camino}
  Let \(\La\subseteq\Z^2\) be an 2-rank integer lattice and \(\vect
  v\in\La\cap\N^2\). There is a lattice path from the origin \(\vect 0\) to
  \(\vect v\) such that no difference between two distinct path nodes
  (except \(\vect v-\vect 0\)) lies in \(\La\) if and only if \(\vect
  v\) is visible in \(\La\) (i.e, the segment joining \(\vect 0\) and
  \(\vect v\) does not contain any other point in \(\La\)) and
  \(\|\vect v\|_1\leq\vol(\La)\).
\end{lemma}

\begin{demo}
  The condition \(\|\vect v\|_1\leq\vol(\La)\) is necessary: the group
  \(\Z^2/\La\) has \(\vol(\La)\) elements, a lattice path with length
  bigger than \(\vol(\La)\) visits at least \(\vol(\La)+2\) points,
  and therefore, at least two pairs of them are congruent modulo
  \(\La\).

  If \(\vect v=(v_1,v_2)\) is not visible in \(\La\), the required
  lattice path cannot exist either. In that case, there is an integer
  \(d\geq 2\) such that \(\frac 1 d\vect v\in\La\). Consider a lattice
  path \(p=(\vect 0=p_0,\ldots,p_{\|\vect v\|_1}=\vect v)\) and its
  subpaths \(p(i,d)\) with origin at \(p_i\) and length \(\frac 1
  d\|\vect v\|_1\), for \(i=0,\ldots,\frac{d-1}d\|\vect v\|_1\). Let
  \(n_i\) be the number of steps of type \(\rightarrow\) that
  \(p(i,d)\) consists of. If, for some index \(i\), we have
  \(n_i=v_1/d\), the difference of the endpoints of \(p(i,d)\) is
  \(\frac 1 d\vect v\in\La\). In other case, since
  \(n_{i+1}-n_i\in\{1,0,-1\}\), we have: \(n_i<v_1/d\), for every
  index \(i\); or \(n_i>v_1/d\), for every index \(i\). Then,
\[v_1=\sum_{i=0}^{d-1}n_{\left(\frac i d\|\vect v\|_1\right)}\neq
v_1.\]

It can be interesting to compare this implication with the so-called
Universal Chord Theorem (see~\citep[p. 15]{Rol_76}). For the other, let
\(\vect v\in\La\cap\N^2\) be a point visible in \(\La\) such that
\(\|\vect v\|_1\leq\vol(\La)\). Draw the segment \(s\) that joins the
origin \(\vect 0\) with \(\vect v\) and consider the diagonals
\(d_i:=\{(x,y)\in\R^2\ |\ x+y=i\}\). We define the path
\(p=(p_0,p_1,\ldots,p_{\|\vect v\|_1})\) by choosing \(p_i\) as the
element of \(d_i\) with integer coordinates that minimizes the
distance to \(s\cap d_i\) (see Figure~\ref{fig_diag}). If the
intersection is equally distant from two integer points, we choose the
one with a bigger first coordinate. It is easy to see that
\(p_{i+1}-p_i\in\{(1,0),(0,1)\}\) and \(p\) is a lattice path from
\(\vect 0\) to \(\vect v\).

\begin{figure}[h!]
\begin{center}
\subfloat[]{\label{fig_diag}\includegraphics[height=4cm]{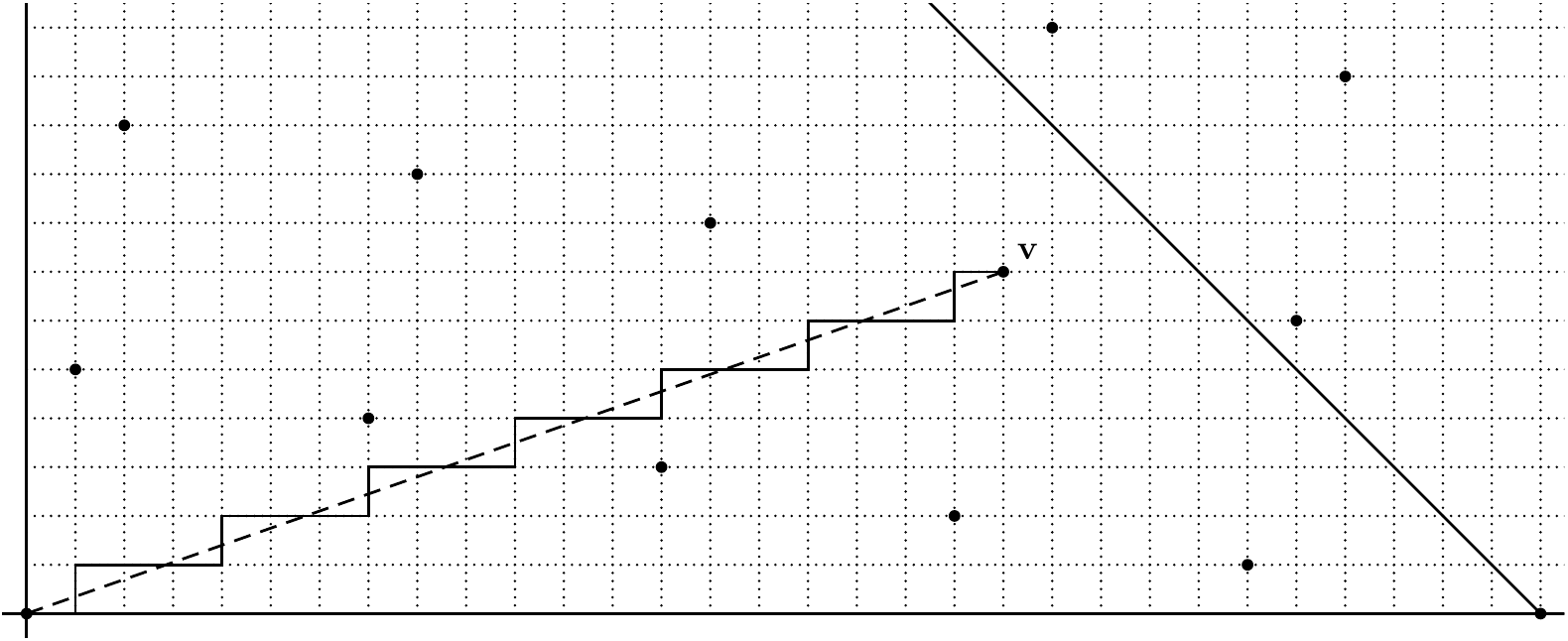}}\hspace{5mm}
\subfloat[]{\label{proj}\includegraphics[height=4cm]{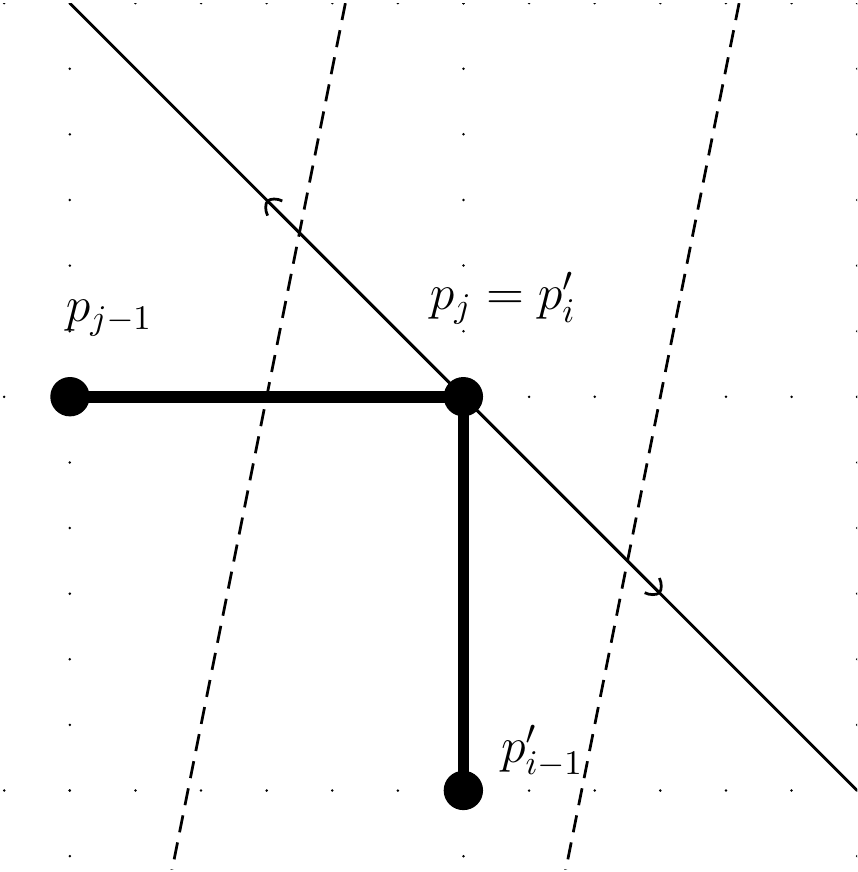}}
\caption{}
\end{center}
\end{figure}

Suppose that \(0\leq i<j\leq\|\vect v\|_1\), \(j-i<\|\vect v\|_1\),
and \(p_j-p_i\in\La\). Note that \(p_j-p_i\) and \(\vect v\) must be
linearly independent, for \(\|p_j-p_i\|_1=j-i<\|\vect v\|_1\) and
\(\vect v\) is visible in \(\La\). Consider the lattice path \(p'\)
``parallel'' to \(p\) with origin at \(p_j-p_i\), i.e.
\(p'_k:=p_k+p_j-p_i\), for \(k=0,\ldots,\|\vect v\|_1\). These paths
collide, for \(p_i'=p_j\). This implies that the image of \(\La\)
through the projection onto \(\R\langle(1,-1)\rangle\) parallel to
\(s\) has two points with distance smaller than \(\sqrt 2\) (a pair of
points distant exactly \(\sqrt 2\) is not enough, as we have set up a
rule for tie breaking).

As \(\vect v\) is visible in \(\La\), there exists \(\vect w\in\Z^2\)
such that \(\La=\Z\langle\vect v,\vect w\rangle\). The projection
described above is defined by the equation:
\[\pi(\vect x)=\vect x-\frac{x_1+x_2}{\|\vect v\|_1}\vect v,\]
so that \(\pi(\La)=\Z\langle\pi(\vect
w)\rangle=\Z\langle\frac{\vol(\La)}{\|\vect v\|_1}(1,-1)\rangle\), and
the distance between two projected lines is at least
\(\frac{\vol(\La)}{\|\vect v\|_1}\sqrt 2\). Therefore the paths
cannot collide and the proposed lattice path satisfies the requirements.
\end{demo}

We close this section with a digression, obtaining, as a corollary of
the previous result, the well-known characterization of circulant
digraphs with two jumps which are Hamiltonian. \citet{Cur_Wit_85} used
similar arguments for the study of Hamiltonian paths.

Consider a connected circulant digraph \(\vec C(n;a,b)\), i.e.
\(\gcd(n,a,b)=1\). A Hamiltonian circuit in this graph corresponds to
a lattice path starting at the origin \(\vect 0\); ending at a certain
\(\vect v\); and such that every element in \(\Z/n\Z\) is the label of
exactly one visited point, except 0, which labels both \(\vect 0\) and
\(\vect v\). In particular, \(\|\vect v\|_1=n\). As the volume of the
circuit lattice of \(\vec C(n;a,b)\) is \(n\), Lemma~\ref{lema_camino}
gives the following characterization of Hamiltonian circulant
digraphs:

\begin{cor}\label{cor_Ham}Let \(n\) be a positive integer and \(a,b\in\Z\) such that
  \(\gcd(a,b,n)=1\). Let \(\La:=\{(u,v)\in\Z^2\ |\ au+bv\in
  n\Z\}\). The circulant digraph \(\vec C(n;a,b)\) is Hamiltonian if
  and only if the intersection \(\La\cap\{(x,y)\in\N^2\ |\ x+y=n\}\)
  has a point visible in \(\La\).
\end{cor}

Note that the intersection of the circuit lattice with the diagonal
\(\{x+y=n\}\) and the first quadrant is \(\{(n-i\frac n d,i\frac n d)\
|\ i=0,\ldots,d\}\), where \(d:=\gcd(b-a,n)\). The point \((n-i\frac n
d,i\frac n d)\) is visible in \(\La\) if and only if the label of
\((d-i,i)\) has order \(n/d\) in \(\Z/n\Z\). Therefore
Corollary~\ref{cor_Ham} is equivalent to the following
characterization (see
\citep{Fio_Yeb_88,Loc_Wit_99,Yan_Bur_Rai_Cel_97}), whose proof dates
back to \citep{Ran_48}.

\begin{prop}Let \(n\) be a positive integer and \(a,b\in\Z\) such that
  \(\gcd(a,b,n)=1\). The circulant digraph \(\vec C(n;a,b)\) is
  Hamiltonian if and only if there exist nonnegative integers \(i,j\) such
  that:
\[i+j=\gcd(b-a,n)=\gcd(ai+bj,n).\]
\end{prop}

\clearpage
\section{Main result}\label{sec_main}

\begin{figure}[h!]
\begin{center}
\subfloat[]{\includegraphics{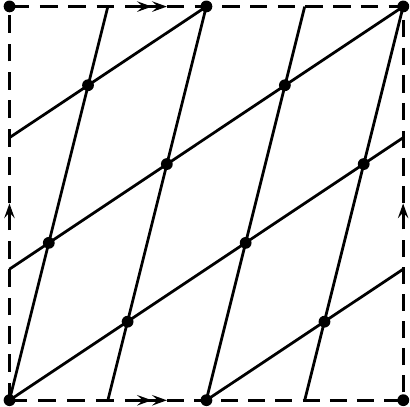}}\hspace{1cm}
\subfloat[]{\label{fig_16}\includegraphics{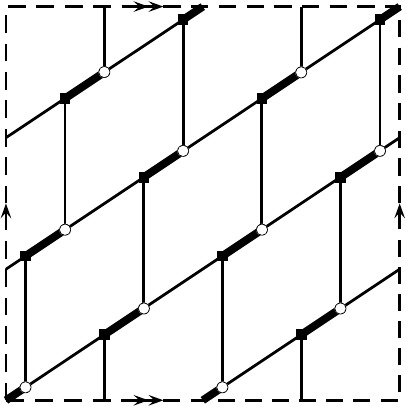}}
\caption{\(\vect v=(3,2),\ \vect w=(1,4)\)}
\end{center}
\end{figure}

As we have seen in Section~\ref{sec_hei}, given a bipartite periodic
graph \(G\) and a matching \(\omega_0\) on its projection \(\hat G\),
any matching \(\omega\in\M(\hat G)\) defines a point
\(h_{\omega_0}(\omega)\in\Z^2\), related to its height
function. Consider the associated Newton polygon \(N_{\omega_0}(G)\):
the convex hull of those points. We prove that any point with integer
coordinates in this polygon is \(h_{\omega_0}(\omega)\), for some
\(\omega\in\M(\hat G)\). This problem can be reduced to visible points
in a triangle with one vertex at the origin and the others at visible
lattice points (named \(\vect v\) and \(\vect w\)). An example of a
graph whose Newton polygon attains this shape is depicted in
Figure~\ref{fig_16}. It is obtained by projecting the segments
\([\vect 0,\vect v]\) and \([\vect 0,\vect w]\) onto the torus and
substituting an \(\omega_0\)-edge for every crossing.
\captionsetup{labelsep=colon}
\captionsetup{labelsep=none}

\begin{figure}[h!]
\begin{center}
\includegraphics[height=4cm]{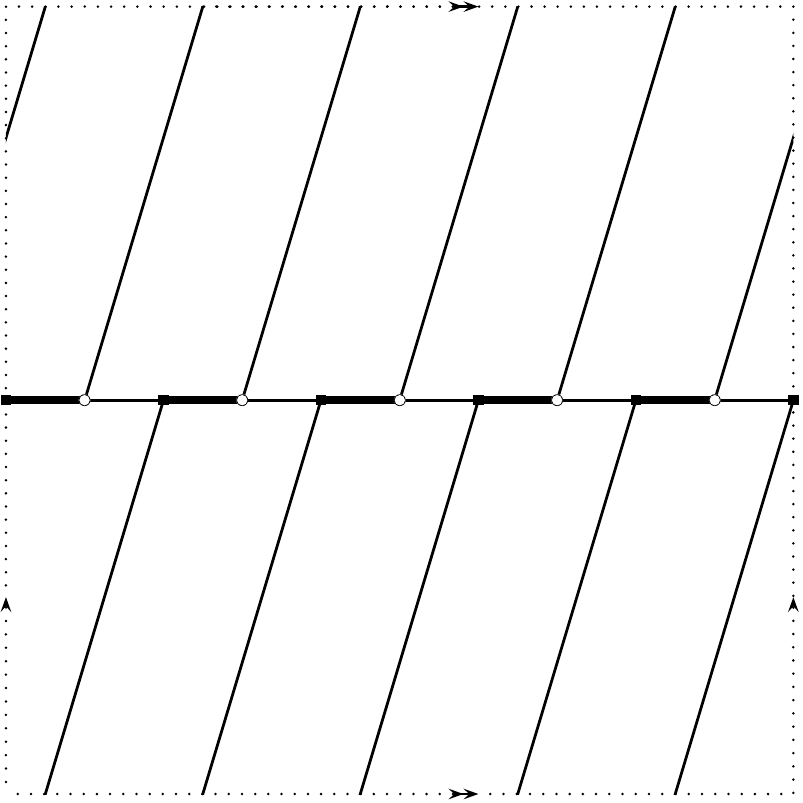}
\caption{A projection of \(B(5,2)\)}\label{est_2}
\end{center}
\end{figure}

For integers \(n,r\) such that \(0\leq r<n\), we consider the
bipartite periodic graph \(B(n,r)\) with black vertices at points
\((i/n,j)\) and white vertices at points \(((2i+1)/2n,j)\), for
integers \(i,j\). Each white vertex \(((2i+1)/2n,j)\) has the three
following neighbours (see Figure~\ref{est_2}): \(\left(\frac i
  n,j\right),\left(\frac{i+1}n,j\right),\left(\frac{i+r}n,j+1\right)\). The
edges of the first type constitute a perfect matching \(\omega_0\),
which we fix as base matching.
\captionsetup{labelsep=colon}
\captionsetup{labelsep=none}

\begin{prop}\label{caso_est}
  Let \(n,r\in\Z\) such that \(0\leq r<n\) and \(\gcd(r,n)=1\), and
  consider the (closed) triangle \(T\) with vertices \((0,0)\),
  \((1,0)\), and \((r,n)\). If \(\vect u\) is a visible lattice point
  in \(T\), there is a matching \(\omega\in\M(\hat B(n,r))\) such that
  \(h_{\omega_0}(\omega)=\vect u\).
\end{prop}
\begin{demo}
  Consider the projected graph \(\hat B(n,r)\) and identify vertex
  \((i/n,j)\) with \(((2i+1)/2n,j)\), obtaining the circulant graph
  \(C(n;1,r)\). Matchings in \(\hat B(n,r)\) correspond with sets of
  disjoint circuits in \(\vec C(n;1,r)\), whose circuit lattice is
  \(\La:=\Z\langle(n,0),(-r,1)\rangle\), and the following bijection
  relates the homology type of a matching to the abelianized of the
  corresponding set of circuits.
\[\begin{array}{ccccc}
  t:&T\cap\Z^2&\longrightarrow&\Lambda\cap\{(x,y)\in\N^2\ |\ x+y\leq
  n\}\\ &\left(\begin{array}{cc}x\\ y\end{array}\right)&\mapsto&\left(\begin{array}{cc}n&-r\\ 0&1\end{array}\right)\left(\begin{array}{cc}x\\ y\end{array}\right)\end{array}\]
Let \(\vect u=(u_1,u_2)\in\Z^2\) be a visible point in
\(T\cap\Z^2\). By Lemma~\ref{lema_camino}, there is a circuit in
\(\vec C(n;1,r)\) whose abelianized is \(t(\vect u)\). If \(\omega\)
is the associated matching, the homology type of \(\omega-\omega_0\)
is \(\vect u\).
\end{demo}

\begin{theo}\label{main_theo}Let \(G\) be a bipartite periodic graph and
  \(\omega_0\in\M(\hat G)\). We have:
\[N_{\omega_0}(G)\cap\Z^2=\{h_{\omega_0}(\omega)\ |\ \omega\in\M(\hat G)\}.\]
\end{theo}

\begin{demo}
  We have to prove that the set on the left-hand side is contained in
  the other, which we denote by \(h(G)\). Let \(\vect u\in
  N_{\omega_0}(G)\cap\Z^2\). It is easy to see that this element is in
  \(h(G)\) if it is a vertex of \(N_{\omega_0}(G)\) or \(\vect u=\vect
  0=h_{\omega_0}(\omega_0)\).

Note that, for any three matchings
\(\omega_1,\omega_2,\omega\in\M(\hat G)\),
\(h_{\omega_1}(\omega)=h_{\omega_1}(\omega_2)+h_{\omega_2}(\omega)\). Therefore
a change in the base matching is equivalent to a translation in the
corresponding set \(h(G)\). Suppose that \(\vect v,\vect w\in
h(G)\). Using Lemma~\ref{lem_divide}, any other point with integer
coordinates in the segment they define lies also in \(h(G)\).

By Caratheodory's theorem, \(\vect u\) lies in the triangle defined by
\(\vect 0\) and two vertices of \(N_{\omega_0}(G)\). Using the
previous argument, we can assume w.l.o.g. that \(\vect u\) is a
visible point in the triangle defined by \(\vect 0\) and two visible
and linearly independent points \(\vect v,\vect w\in h(G)\). By
Lemma~\ref{lem_divide_bis}, there exist \(\omega_1,\omega_2\in\M(\hat
G)\) with respective homology types \(\vect v\) and \(\vect w\), and
such that each of the transition graphs \(\omega_1-\omega_0\) and
\(\omega_2-\omega_0\) consists of a single circuit.

We know (see the discussion in Subsection~\ref{sec_top}) that there is
a self-homeomorphism \(f\) on the torus whose associated automorphism
on the first homology group satisfies:
\[\pi_1(f)(\vect v|\vect w)=\left(\begin{array}{cc}
    1&r\\ 0&n\end{array}\right),\] where \(0\leq r<n:=|\det(\vect v|\vect w)|\)
and \(\gcd(r,n)=1\). Using Lemma~\ref{lem_b}, another torus
self-homeomorphism transforms the transition cycle
\(\omega_1-\omega_0\) into the meridian \([0,1]\times\{0\}\). The set of height change vectors of the transformed
graph \(\hat G_1\) is \(\pi_1(f)h(G)\), so we need just prove that
any visible point in the triangle defined by \((0,0),(1,0),\) and
\((r,n)\) is in \(h(G_1)\).

Now, we transform \(\hat G_1\) into another graph \(\hat H\) such
that \(h(H)\subseteq h(G_1)\). Consider the intersections of the
meridian (image of \(\omega_1-\omega_0\)) and the image of the circuit
\(\omega_2-\omega_0\). Such an intersection cannot be limited to a single
point.
\begin{figure}[h!]
\begin{center}
\subfloat[]{\label{fig_incident}\includegraphics{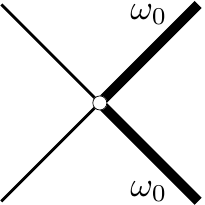}}\hspace{15mm}\vline\hspace{15mm}
\subfloat[]{\label{fig_cap}\includegraphics{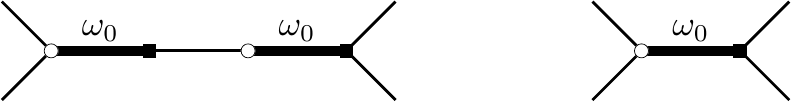}}
\caption{}
\end{center}
\end{figure}
If this was the case, that point should be a vertex incident to two
different edges of \(\omega_0\), which would not be a perfect matching
(see Figure~\ref{fig_incident}). Therefore each connected component
of the intersection consists of an odd number of edges, which can be
reduced to one, as depicted in Figure~\ref{fig_cap}.

Suppose that in \(G_1\) there are two consecutive (as walking through
the meridian) crossings between both considered circuits with the same
\(Y\)-coordinate.
\begin{figure}[h!]
\begin{center}
\includegraphics{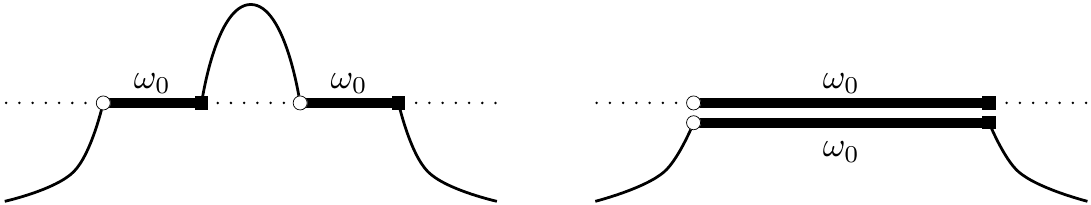}
\caption{}\label{plancha}
\end{center}
\end{figure}
We can remove (see Figure~\ref{plancha}) all those crossings (which
happen in finite number, for the projected graph is finite), obtaining
a graph \(H\) homeomorphic (as a planar embedding) to
\(B(n,r)\). Indeed, after the removal process, there are exactly \(n\)
crossings in \(\hat H\), and the edges of the second circuit define a
permutation of them. As there is no self-intersection, this
permutation must be of the type \(x\mapsto x+s\). We have
\(\gcd(s,n)=1\) (\(s=r\), indeed), for \(\omega_2-\omega_0\) consists
of a single circuit. Therefore Proposition~\ref{caso_est} applies.
\end{demo}

\noindent{\bf Acknowledgements.} The author is supported by the
Spanish \emph{Ministerio de Educación, Cultura y Deporte}, through the
``Programa Nacional de Movilidad de Recursos Humanos del Plan Nacional
de I+D+i 2008-2011'' and the Spanish \emph{Ministerio de Economía y
Competitividad}, through grant MTM2011-24678.

\bibliography{FSK}

\end{document}